# AIDE: Fast and Communication Efficient Distributed Optimization

Sashank J. Reddi  
Carnegie Mellon University

Jakub Konečný  
University of Edinburgh

Peter Richtárik  
University of Edinburgh

Barnabás Póczós  
Carnegie Mellon University

Alex Smola  
Carnegie Mellon University

**Abstract**

In this paper, we present two new communication-efficient methods for distributed minimization of an average of functions. The first algorithm is an inexact variant of the DANE algorithm [20] that allows any local algorithm to return an approximate solution to a local subproblem. We show that such a strategy does not affect the theoretical guarantees of DANE significantly. In fact, our approach can be viewed as a robustification strategy since the method is substantially better behaved than DANE on data partition arising in practice. It is well known that DANE algorithm does not match the communication complexity lower bounds. To bridge this gap, we propose an accelerated variant of the first method, called AIDE, that not only matches the communication lower bounds but can also be implemented using a purely first-order oracle. Our empirical results show that AIDE is superior to other communication efficient algorithms in settings that naturally arise in machine learning applications.

## 1  Introduction

With the advent of large scale datasets, distributed machine learning has garnered significant attention recently. For example, large scale distributed machine learning systems such as the *Parameter server* [9], *GraphLab* [22] and *TensorFlow* [1] work with datasets sizes in the order of hundreds of terabytes. When dealing with datasets of such scale in distributed systems, computational and communication workloads need to be designed carefully. This in turn places primary importance to computational and communication resource constraints on the algorithms used in these machine learning systems.

In this paper, we study the problem of distributed optimization for solving *empirical risk minimization* problems, where one seeks to minimize average of $N$ functions $\min_{w \in \mathbb{R}^d} f(w) := \frac{1}{N} \sum_{i=1}^N f_i(w)$. Many problems in machine learning, such as logistic regression or deep learning, fall into this setting. Formally, we assume a distributed system consisting of $K$ machines, where machine $k$ has access to a (nonempty) subset of the data indexed by $\mathcal{P}_k \subset [N] := \{1, \ldots, N\}$. We assume that $\{\mathcal{P}_k\}$ forms a partition of $[N]$. Our goal can then thus be reformulated as

$$\min_{w \in \mathbb{R}^d} f(w) = \frac{1}{K} \sum_{k=1}^K F_k(w), \text{ where } F_k(w) := \frac{K}{N} \sum_{i \in \mathcal{P}_k} f_i(w). \qquad (1)$$

We assume that (1) has an optimal solution. By $\hat{w}$ we denote an arbitrary optimal solution: $\hat{w} \in \arg\min_w f(w)$. The fundamental constraint is that machine $k$ can directly access only (the data describing) function $F_k(w)$, i.e., functions $f_i(w)$ for $i \in \mathcal{P}_k$. Typically, $f_i(w) = \phi_i(w^\top x_i)$, where $x_i$ is the $i^{th}$ data point, and $\phi_i(\cdot)$ is a loss function (dependence on $i$ indicates potential dependence on a label). The function $f_i$ is assumed to be continuous but *not* necessarily convex — the structure also encompasses problems in deep learning.

The basic benchmark algorithm for solving the optimization problem is gradient descent (GD). Each iteration of GD performs the update step $w \leftarrow w - \frac{h}{K} \sum_{k=1}^{K} \nabla F_k(w)$, where $h > 0$ is a stepsize parameter. In a distributed system, this amounts to calculation of the gradient $\nabla F_k(w)$ by each machine $k$ and then sending an update to a central node to compute the full update direction $\frac{1}{K} \sum_{k=1}^{K} \nabla F_k(w)$. Although such an update sequence is simple, it is often communication intensive due to the communication involved at each iteration of the algorithm. Furthermore, because of the communication delay, computational resources at each machine are often under-utilized. While the latter problem can be addressed by using asynchronous and stochastic variants of gradient descent [15, 16], these algorithms still incur high communication costs at each iteration.

Recently, there have been considerable theoretical advances in our understanding of the communication overheads in solving (1). In particular, three different methods address this issue by designing a procedure that uses significant amount of computation locally, between rounds of communication. These methods are DANE [20], DISCO [23] and COCOA+ [11, 5, 12, 21]. Both DANE and DISCO show that when the functions $\{F_k\}_{k=1}^{K}$ in (1) are similar (see Definition 1) or exhibit special structure, one could obtain the optimal solution in considerably fewer rounds of communication.

At the $t^{\text{th}}$ iteration of DANE, the following general subproblem has to be solved *exactly* by machine $k$, with $\eta \geq 0$ and $\mu \geq 0$ being parameters of the method:

$$\min_{w \in \mathbb{R}^d} F_k(w) - \langle \nabla F_k(w^{t-1}) - \eta \nabla f(w^{t-1}), w \rangle + \frac{\mu}{2} \|w - w^{t-1}\|^2,$$

where $\langle \cdot, \cdot \rangle$ is the standard inner product and $\|\cdot\| := \langle \cdot, \cdot \rangle^{1/2}$ is the standard Euclidean norm. The gradient $\nabla f(w^{t-1})$ is computed in a distributed fashion at the start, followed by aggregation of individual updates. DISCO uses an *inexact* damped Newton method for solving the problem (1) [23]; and hence is a *second* order method. Thus, both DANE and DISCO require strong oracle access — either in the form of oracle to solve another optimization problem exactly, or second order information. This is in sharp contrast with COCOA+, where (similarly to DANE) a general local problem is formulated, but is needed to be solved only approximately. This enables use of essentially any algorithm to be used locally for a number of iterations, making the method much more versatile.

A natural idea is to solve DANE subproblem approximately using a first-order method. However, from the point of analysis, this leads to an additional error. Hence, it is not clear whether such method enjoys any theoretical guarantees, as such a generalization is non-trivial. Based on this intuition, we develop an inexact version of DANE (refereed to as INEXACTDANE) and provide its convergence analysis. We demonstrate that our proposed approach is significantly more robust to 'bad' data partitioning, and also highlight a connection to a distributed version of SVRG algorithm [6, 8]; thus, yielding partial convergence guarantees for a method that has been observed to perform well in practice [7], but has not been successfully analyzed.

While INEXACTDANE is practically appealing in comparison to DISCO due to its simplicity and possible reliance only on first-order information, DISCO is *theoretically* superior to INEXACTDANE in terms of communication complexity, as it nearly matches the lower bounds derived in [2]. To address this issue, we build upon INEXACTDANE and propose an approach, referred to as AIDE, that matches (up to logarithmic factors) these communication complexity lower bounds and can be implemented using a first-order oracle. To our knowledge, ours is the first work that can be implemented using a first-order oracle, but at the same time achieves near optimal communication complexity for the setting considered in this paper.

## 1.1 Related Work

The literature on distributed optimization is vast and hence, we restrict our attention to the works most relevant to our paper. As mentioned earlier, the most straightforward approach for solving the distributed optimization is using gradient descent or its accelerated version. When the function $f$ is $L$-smooth and $\lambda$-strongly convex, the communication complexity of accelerated gradient descent is $O(\sqrt{L/\lambda} \log(1/\epsilon))$ for achieving $\epsilon$-suboptimal solution [13]. Another possible approach, often referred to as "one shot averaging", is solve the local optimization problem parallelly at each machine and then compute average of the solutions [25, 24]. However, it can be shown that "one shot averaging" can produce significantly worse solution than the optimal solution $\hat{w}$.

ADMM (alternating direction method of multipliers) is another popular approach for solving distribution optimization problems. Under certain conditions, ADMM is shown to achieve communication complexity of $O(\sqrt{L/\lambda}\log(1/\epsilon))$ for $L$-smooth and $\lambda$-strongly convex function. More recently, the above mentioned DANE, DISCO and COCOA+ algorithms have been proposed to tackle the problem of reducing the communication complexity in solving problems of form (1). We defer a detailed comparison of our algorithms with DANE and COCOA+ to Section 7 and Appendix E. The lower bounds on the communication complexity for solving problems of form (1) have been obtained in [2] (see Section 8 for more details). We refer the readers to [20, 2] for a more comprehensive coverage of the related work.

## 2  Algorithm: INEXACTDANE

The DANE algorithm [20] proceeds as follows. At iteration $t$, node $k$ *exactly* solves the subproblem $\hat{w}_k^t = \arg\min_w g_{k,\mu}^t(w)$, where

$$g_{k,\mu}^t(w) := F_k(w) - \langle \nabla F_k(w^{t-1}) - \eta \nabla f(w^{t-1}), w\rangle + \tfrac{\mu}{2}\|w - w^{t-1}\|^2, \tag{2}$$

and $\mu, \eta \geq 0$ are parameters. After this, the method computes $w^t = \sum_{k=1}^K \hat{w}_k^t/K$, which involves communication, and the process is repeated. Our first method, INEXACTDANE, modifies this algorithm by allowing the subproblems to be solved inexactly. In particular, let $w_k^t$ denote the point obtained by minimizing $g_{k,\mu}^t$ only approximately (for example, using Quartz [14] or SVRG [6, 8]). The pseudocode of this method is formalized in Algorithm 1. The parameter $\gamma$ refers to the *level of inexactness* allowed (small $\gamma$ means higher accuracy). Note that each iteration of INEXACTDANE can be solved in an embarrassingly parallel fashion.

The communication complexity of an algorithm is defined as the number of communication rounds required by the algorithm to reach a solution $w^t$ satisfying specific convergence criteria such as $f(w^t) - f(\hat{w}) \leq \epsilon$. In each communication round, the machines can only send information linear in size of the dimension $d$. This is also the notion of communication round used in [2].

We start with the analysis of *quadratic case* (Section 3), after which we deal with strongly convex, weakly convex and nonconvex cases (Section 4).

## 3  Analysis of INEXACTDANE: Quadratic Case

In this section we assume that $F_k(w) = \tfrac{1}{2}w^\top H_k w + l_k^\top w$, where $H_k \in \mathbb{R}^{d\times d}$ and $l_k \in \mathbb{R}^d$ (this is the case when the functions $\{f_i\}_{i=1}^N$ in (1) are quadratic). We write $H := \tfrac{1}{K}\sum_{k=1}^K H_k$ and $l := \tfrac{1}{K}\sum_{k=1}^K l_k$. For a square matrix $A$, by $\|A\|$ we denote its spectral norm. Our first result plays a central role in the analysis of quadratic functions. Proofs are provided in the Appendix.

**Theorem 1.** *Assume $H \succ 0$. Let $\mu \geq 0$ and define $\tilde{H}^{-1} := \tfrac{1}{K}\sum_{k=1}^K (H_k + \mu I)^{-1}$. Further, choose $0 \leq \gamma < 1$, $\eta > 0$ and define*

$$\rho := \|\eta\tilde{H}^{-1}H - I\| + \tfrac{\eta\gamma}{K}\sum_{k=1}^K \|(H_k + \mu I)^{-1}H\|.$$

*Let $\{w^t\}$ be the iterates of INEXACTDANE (Algorithm 1 applied with Option I). Then for all $t \geq 1$ we have $\|w^t - \hat{w}\| \leq \rho \|w^{t-1} - \hat{w}\|$.*

Suppose, $0 < \lambda \preceq H_k \preceq L$ for all $k \in [K]$, then with $\eta = 1$, sufficiently large $\mu$ and small $\gamma$, one can ensure that $\rho < 1$. We are interested in the setting where the functions $\{F_k\}$ are "similar". Following [2], we shall measure similarity via the notion of $\delta$-*relatedness*, defined next.

**Definition 1.** *We say that quadratic functions $\{F_k\}$ are $\delta$-related, if for all $k, k' \in [K]$ we have*

$$\|H_k - H_{k'}\| \leq \delta, \qquad \|l_k - l_{k'}\| \leq \delta.$$

**Algorithm 1:** INEXACTDANE($f, w^0, s, \gamma, \mu$)

**Data:** $f(w) = \frac{1}{K}\sum_{k=1}^{K} F_k(w)$, initial point $w^0 \in \mathbb{R}^d$, inexactness parameter $0 \leq \gamma < 1$
**for** $t = 1$ **to** $s$ **do**
    **for** $k = 1$ **to** $K$ **do in parallel**
        Find an approximate solution $w_k^t \approx \hat{w}_k^t := \arg\min_{w \in \mathbb{R}^d} g_{k,\mu}^t(w)$
        **Option I:** $\|w_k^t - \hat{w}_k^t\| \leq \gamma \|w^{t-1} - \hat{w}_k^t\|$
        **Option II:** $\|\nabla g_{k,\mu}^t(w_k^t)\| \leq \gamma \|\nabla g_{k,\mu}^t(w^{t-1})\|$ and $\|w_k^t - \hat{w}_k^t\| \leq \gamma \|w^{t-1} - \hat{w}_k^t\|$
    **end**
    $w^t = \frac{1}{K}\sum_{k=1}^{K} w_k^t$
**end**
**return** $w^t$

The following theorem specifies the convergence rate in the case of $\delta$-related objectives. We are particularly interested in the setting where $\delta$ is small, which is the case that should allow stronger communication complexity guarantees [20, 2].

**Theorem 2.** *Let $0 < \lambda \leq L$ be such that $\lambda I \preceq H \preceq LI$, and assume that $\{F_k\}$ are $\delta$-related, with $\delta \geq 0$. Choose $\mu = \max\{0, (8\delta^2/\lambda) - \lambda\}$, $\eta = 1$ and let*

$$\gamma = \begin{cases} 1/8 & \text{if } 2\sqrt{2}\delta \leq \lambda \\ \lambda^2/(192\delta^2) & \text{otherwise} \end{cases} \qquad \tilde{\rho} := \begin{cases} 2/3 & \text{if } 2\sqrt{2}\delta \leq \lambda \\ 1 - (\lambda^2/24\delta^2) & \text{otherwise} \end{cases}$$

*The iterates $\{w^t\}$ of* INEXACTDANE *(Algorithm 1 with Option I) satisfy $\|w^t - \hat{w}\| \leq \tilde{\rho}\|w^{t-1} - \hat{w}\|$.*

It is interesting to note that solving the local subproblems beyond the accuracy $\gamma$ stated in the above result does *not* lead to a better overall complexity result. That is, the required accuracy at the nodes of the distributed system should not be perfect, but should instead depend in $\delta$ and $\lambda$. We have the following corollary, translating the result into a bound on the number of iterations (i.e., number of communication rounds) sufficient to obtain an $\epsilon$-solution.

**Corollary 1.** *Let the assumptions of Theorem 2 be satisfied, and pick $0 < \epsilon < L\|w^0 - \hat{w}\|^2$. If*

$$t = \tilde{O}\left(\frac{\delta^2}{\lambda^2}\log\left(\frac{L\|w^0 - \hat{w}\|^2}{\epsilon}\right)\right),$$

*then $f(w^t) \leq f(\hat{w}) + \epsilon$.*

The $\tilde{O}(\cdot)$ notation hides few logarithmic factors. In the situation when for all $k$, $H_k$ arises as the average of Hessians of $n$ i.i.d. quadratic functions with eigenvalues upper-bounded by $L$, it can be shown that $\delta = O(L/\sqrt{n})$, hence $t = \tilde{O}((L/\lambda)^2/n \log(1/\epsilon))$ (see Corollary 4 in the Appendix). This arises for instance in a linear regression setting where the samples at each machine are i.i.d — a scenario typically assumed in distributed machine learning systems. Note that, in this case, the communication complexity *decreases* as $n$ increases.

## 4 Analysis of INEXACTDANE: General Case

In this section, we present the results for general strongly convex case, weakly convex case, and non-convex case. Proofs are provided in the Appendix.

A function $\psi : \mathbb{R}^d \to \mathbb{R}$ is called *L-smooth* if it is differentiable, and if for all $x, y \in \mathbb{R}^d$ we have $\|\nabla\psi(x) - \nabla\psi(y)\| \leq L\|x - y\|$. It is called $\lambda$-*strongly convex* if $\psi(x) \geq \psi(y) + \langle\nabla\psi(y), x - y\rangle + \frac{\lambda}{2}\|x - y\|^2$, for all $x, y \in \mathbb{R}^d$ (if $\psi$ is twice differentiable, this is equivalent to requiring that $\nabla^2\psi(w) \succeq \lambda I$ for all $w \in \mathbb{R}^d$).

If $\lambda > 0$, we say that $\psi$ is strongly convex. In this case, $\kappa = L/\lambda$ denotes the condition number. We say that $\psi$ is *weakly convex* if it is 0-strongly convex.

## 4.1 Strongly convex case

Our analysis of the strongly convex case follows along the lines of [20] and incorporates the inexactness in solving the subproblem in (2) at each iteration.

**Theorem 3.** *Assume that for all $k \in [K]$, $F_k$ is L-smooth and $\lambda$-strongly convex, with $\lambda > 0$. Let*

$$\tilde{\rho} := \left[ \frac{(1-\gamma)^2}{\eta(L+\mu)} - \frac{2L}{(\lambda+\mu)^2} - \frac{2\gamma(L+\mu)}{\eta(\lambda+\mu)^2} \right] \eta^2 \lambda. \tag{3}$$

*Suppose $\eta > 0$, $0 \leq \gamma < 1$ and $\mu > 0$ are chosen such that $0 < \tilde{\rho} < 1$. Then for the iterates of* INEXACTDANE *(Algorithm 1 with Option II) we have:* $f(w^t) - f(\hat{w}) \leq (1 - \tilde{\rho}) \left( f(w^{t-1}) - f(\hat{w}) \right)$.

By setting $\gamma = 0$, we roughly recover the corresponding result in [20] covering the exact case. Note that $\gamma$ only has a weak effect on the convergence rate. For instance, with $\gamma = 1/8$, $\eta = 1$ and $\mu = 6L - \lambda$, we require $O(\frac{L}{\lambda} \log(1/\epsilon))$ iterations to achieve a solution $w^t$ such that $f(w^t) - f(\hat{w}) \leq \epsilon$.

## 4.2 Weakly convex case

We analyze the weakly convex case by using a perturbation argument. This allows us to use the analysis of strongly convex case for proving the convergence analysis. In particular, we consider the following function $f_\epsilon(w) = f(w) + \frac{\epsilon}{2} \|w - w^0\|^2$, which is essentially a perturbation of the function $f$. First, note that if $f$ is $L$-smooth then $f_\epsilon$ is $(L + \epsilon)$-smooth and $\epsilon$-strongly convex. Applying INEXACTDANE (Algorithm 1) on the function $f_\epsilon$ and using Theorem 3, we get the following:

$$f_\epsilon(w^t) - \min_w f_\epsilon(w) \leq (1 - \rho_\epsilon)^s \left[ f_\epsilon(w^0) - \min_w f_\epsilon(w) \right] \tag{4}$$

where $\rho_\epsilon$ is defined as in (3) with $\lambda$ and $L$ replaced by $\epsilon$ and $L + \epsilon$ respectively. Using the above relationship, the corollary below follows immediately.

**Corollary 2.** *Suppose the function $f$ is weakly convex. Then the iterates of* INEXACTDANE *in Algorithm 1 (Option II) with $\eta = 1$, $\gamma = 1/8$ and $\mu = 6L + 5\epsilon$ applied to $f_\epsilon$ after $s = \tilde{O}(L \log(1/\epsilon)/\epsilon)$ iterations satisfy $f(w^t) \leq f(\hat{w}) + O(\epsilon)$.*

The result essentially shows sublinear convergence rate (up to logarithmic factor) of INEXACTDANE for weakly convex functions and weak dependence on the inexact parameter $\gamma$.

## 4.3 Nonconvex case

Finally, we look at the case where function $f$ can be nonconvex. The key observation is that by choosing large enough $\mu$ one can have a strongly convex $g_{k,\mu}^t$ in (2). The existence of such a $\mu$ is due to the Lipschitz continuous nature of gradient of $f$. This allows us to solve the subproblem efficiently and makes it amenable to analysis. For the purpose of theoretical analysis, in the nonconvex case, we select $w^t$ arbitrarily from $\{w_k^t\}_{k=1}^K$ instead of averaging as in Algorithm 1.

**Theorem 4.** *Assume that for all $k \in [K]$, $F_k$ is L-smooth. Suppose the iterates $w_k^t$ of Algorithm 1 (with Option II) for some $0 \leq \gamma < 1$. Furthermore, let*

$$\theta := \left[ \frac{(1-\gamma)^2}{\eta(L+\mu)} - \frac{2L}{(\mu-L)^2} - \frac{2\gamma(L+\mu)}{\eta(\mu-L)^2} \right] \eta^2,$$

*where $\eta > 0$, $0 \leq \gamma < 1$ and $\mu > L$ are such that $\theta > 0$. Then, we have*

$$\min_{0 \leq t' \leq t-1} \|\nabla f(w^{t'})\|^2 \leq \frac{f(w^0) - f(\hat{w})}{\theta t}.$$

---
**Algorithm 2:** AIDE($f, w^0, \lambda, \tau, s, \gamma, \mu, \epsilon$)

**Data:** $f(w) = \frac{1}{K}\sum_{k=1}^{K} F_k(w)$, Initial point $y^0 = w^0 \in \mathbb{R}^d$, INEXACTDANE iterations $s$, inexactness parameter $0 \leq \gamma < 1, \tau \geq 0$

Let $q = \lambda/(\lambda + \tau)$
**while** $f(w^{t-1}) - f(\hat{w}) \leq \epsilon$ **do**
    Define $f^t(w) := \frac{1}{K}\sum_{k=1}^{K}(F_k(w) + \frac{\tau}{2}\|w - y^{t-1}\|^2)$
    $w^t$ = INEXACTDANE($f^t, w^{t-1}, s, \gamma, \mu$)
    Find $\zeta_t \in (0,1)$ such that $\zeta_t^2 = (1-\zeta_t)\zeta_{t-1}^2 + q\zeta_t$
    Compute $y^t = w^t + \beta_t(w^t - w^{t-1})$ where $\beta_t = \frac{\zeta_{t-1}(1-\zeta_{t-1})}{\zeta_{t-1}^2 + \zeta_t}$
**end**
**return** $w^t$
---

With $\eta = 1$, one could choose $\mu = 10L$ and constant inexactness of $\gamma = 1/8$ to get $\theta \geq 1/100L$. Again, similar to the strongly convex case, we see a very weak effect of $\gamma$ on the convergence rate. In other words, with constant approximation at each local node, we can obtain $O(1/t)$ convergence rate to a stationary point for nonconvex functions.

## 5 Accelerated Distributed Optimization

In this section, we study an accelerated version of INEXACTDANE algorithm. While INEXACTDANE algorithm reduces the communication complexity in distributed settings, it is known that it is *not* optimal i.e., there is a gap between upper bounds of INEXACTDANE and the lower bounds for the communication complexity proved in [2]. To this end, we propose an accelerated variant of INEXACTDANE algorithm and prove that it matches the lower bounds (upto logarithmic factors) in specific settings. We refer to the Accelerated Inexact DanE as "AIDE" ( Algorithm 2). The key idea is to apply the generic acceleration scheme—catalyst—in [10] to INEXACTDANE. We show that by a careful selection of $\tau$ in Algorithm 2, one can achieve optimal communication complexity in interesting and important settings.

### 5.1 Quadratic case

Here we show that by appropriate selection of $\tau$, for quadratic case, one can achieve faster convergence rates that match the lower bounds presented in [2].

**Theorem 5.** *Assume that $\lambda I \preceq H \preceq LI$ and further assume that $\{F_k\}$ are $\delta$-related. Then the iterates of AIDE (Algorithm 2) with $\eta = 1$, $\mu = 0$, $\tau = \max\{0, 2\sqrt{2}\delta - \lambda\}$, $s = \tilde{O}(1)$ and $\gamma = 1/8 - (\delta^2/2(\tau+\lambda)^2)$ after $t = \tilde{O}(\sqrt{\delta/\lambda}\log(1/\epsilon))$ satisfy $f(w^t) \leq f(\hat{w}) + \epsilon$ and the total number of iterations in INEXACTDANE (Option I) is $\tilde{O}(\sqrt{\delta/\lambda}\log(1/\epsilon))$.*

Note that the communication cost of AIDE is proportional to the number of iterations of INEXACTDANE executed as part of AIDE algorithm. Hence, the total communication complexity of AIDE is $\tilde{O}(\sqrt{\delta/\lambda}\log(1/\epsilon))$. It is interesting to note that it is enough to solve each subproblem in INEXACTDANE up to a constant accuracy of $\gamma = 1/8$ in order to achieve the convergence rate stated above. In other words, by using AIDE, one need not solve the subproblems with high accuracy and still achieve faster convergence than INEXACTDANE.

### 5.2 Convex case

For the general strongly convex case, we prove the following key result concerning the number of iterations of AIDE and INEXACTDANE used as part of AIDE.

**Theorem 6.** *Assume that for all $k \in [K]$, the function $F_k$ is $L$-smooth and $\lambda$-strongly convex. Then the iterates of* AIDE *with parameters $\lambda$, $\eta = 1$, $\mu = 12L$, $\gamma = 1/8$, $s = \tilde{O}(1)$ and $\tau = L - \lambda$ after $t = \tilde{O}(\sqrt{L/\lambda} \log(1/\epsilon))$ iterations satisfy $f(w^t) \leq f(\hat{w}) + \epsilon$ and the total number of iterations in* INEXACTDANE *(Option II) is $\tilde{O}(\sqrt{L/\lambda} \log(1/\epsilon))$.*

The upper bound on the communication complexity matches the lower bounds proved in [2] for unrelated strongly convex functions. It is an interesting open problem to obtain better communication complexity for a subclass of strongly convex functions other than the quadratic case explored above.

As in the case of INEXACTDANE, we apply AIDE to a perturbed problem when $f$ is weakly convex. In this manner, we invoke catalyst (acceleration) for strongly convex functions. Recall the function $f_\epsilon$ is defined as $f_\epsilon(w) = f(w) + \frac{\epsilon}{2} \|w - w^0\|^2$. Using AIDE on $f_\epsilon$, Theorem 6 shows that one can achieve $\epsilon$-accuracy on function $f_\epsilon$ in $t = \tilde{O}(\sqrt{(L+\epsilon)/\epsilon} \log(1/\epsilon))$ iterations. This follows from the fact that $f_\epsilon$ is $L + \epsilon$-smooth and $\epsilon$-strongly convex. Again, using similar argument as in Equation (17), we have the following corollary.

**Corollary 3.** *Suppose the function $f$ is weakly convex. Then iterates of* AIDE *in Algorithm 2 applied on $f_\epsilon$ with $\lambda = \epsilon$, $\eta = 1$, $\gamma = 1/8$ and $\mu = 12(L+\epsilon)$, $s = \tilde{O}(1)$ and $\tau = L$ after $t = \tilde{O}(\sqrt{L/\epsilon} \log(1/\epsilon))$ iterations satisfy $f(w^t) \leq f(\hat{w}) + O(\epsilon)$.*

## 6 Connection to a Practical Distributed Version of SVRG

In this section, we discuss the connection between SGD, INEXACTDANE and a version of distributed SVRG that was observed to perform well in practice [7], but has not yet been successfully analyzed. Algorithm 1 uses an approximate solution to the minimization of function $g_{k,\mu}^t(w)$ at each iteration; however, it does not specify the algorithm (local solver) used to obtain it. While one could use any of recent fast algorithms for finite sums [17, 4, 19, 3] with comparable results, we highlight the consequences of applying the SVRG algorithm [6, 8] as a local solver in INEXACTDANE.

The SVRG algorithm for minimizing average of functions $f(w) = \frac{1}{N} \sum_{i=1}^{N} f_i(w)$ runs in two nested loops. In the $t^{\text{th}}$ outer iteration, SVRG computes the gradient of the entire function $\nabla f(w^{t-1})$ followed by inner loop, where an update direction at $w$ is computed as $\nabla f_i(w) - \nabla f_i(w^{t-1}) + \nabla f(w^{t-1})$ for a randomly chosen index $i \in [N]$. For its distributed version, the inner loop is executed parallelly over all the machines but the index $i$ for machine $k$ is sampled only from $\mathcal{P}_k$ (the data available locally) instead of $[N]$. For completeness, the pseudocode is provided in Appendix D.

We can obtain an identical sequence[1] of iterates $w^t$ by looking at a variation of INEXACTDANE algorithm. It is easy to verify that if we apply SVRG to the INEXACTDANE subproblem $g_{k,\mu}^t$ in (2) with $\mu = 0$ and $\eta = 1$, running SVRG for a single outer iteration, the resulting procedure is equivalent to running the Algorithm 3 (in the Appendix). Pointing out this connection gives partial justification for the performance of Algorithm 3. However, a direct analysis for Algorithm 3 still remains an open question as the above reasoning applies only when $\mu = 0$ and $\eta = 1$. Note that our results apply to this special setting only under specific conditions; for example when quadratic functions $\{F_k\}$ that are $\delta$-related with sufficiently small $\delta$.

## 7 Experiments

In this section, we demonstrate the empirical performance of the INEXACTDANE and AIDE algorithms on binary classification task using different loss functions. We compare the performance of our algorithms to that of COCOA+, a popular distributed algorithm [11]. COCOA+ is particularly relevant to our setting since, similar to INEXACTDANE and AIDE, it provides flexibility in dealing with communication/computation tradeoff by choosing the local computation effort spent between rounds of communication. We use SVRG as the local solver for INEXACTDANE and AIDE, and SDCA as the local solver for COCOA+ in all our experiments. Unless stated otherwise, at each iteration of the algorithms, we run the corresponding local solvers for a single pass over data available locally. Note that for all the algorithms considered here, the iteration complexity is proportional to their communication complexity.

---
[1] Subject to the same sequence on sampled indices $i \in \mathcal{P}_k$.

For our experiments, we use standard binary classification datasets rcv1, covtype, real-sim and url[2]. As part of preprocessing, we normalize the data and add a bias term. In our experiments, the data is randomly partitioned across individual nodes; thus, mimicking the i.i.d data distribution. We are interested in examining the effect of varying amount of local computation and number of nodes on the performance of the algorithms. For our results, we present the best performance obtained by selected from a range of stepsize parameters for SVRG and aggregation parameters for COCOA+. For COCOA+, this amounts to searching for the optimal choice of parameter $\sigma'$. In the plots, we use DANE as the label for INEXACTDANE.

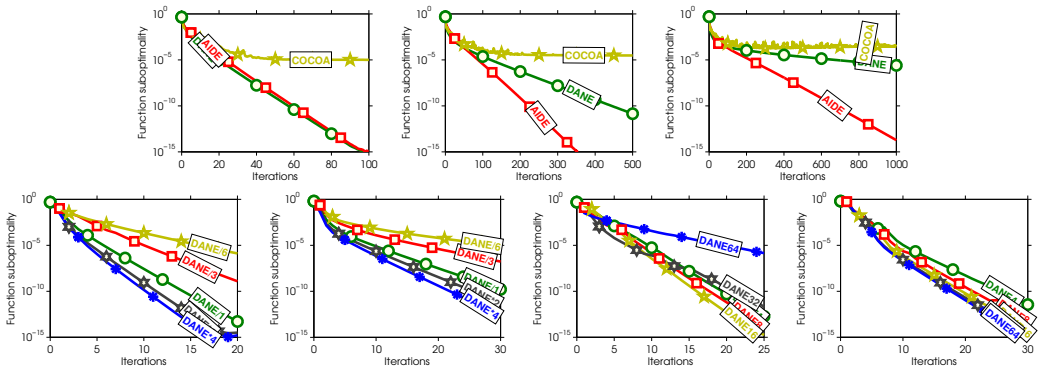

Figure 1: Top: rcv1 dataset; Smoothed hinge loss; $\lambda$ set to $1/(cN)$, for $c \in \{1, 10, 100\}$. Bottom: Logistic loss; *left-hand two:* Varying number of local data passes per iteration. rcv1 and url datasets. *right-hand two:* Varying number of nodes; rcv1 and url datasets.

In the top row of Figure 1, we compare INEXACTDANE, AIDE and COCOA+ for classification task on the rcv1 dataset with smoothed hinge loss on 8 nodes. The three plots are with regularization parameter $\lambda$ set to $1/(cN)$, for $c \in \{1, 10, 100\}$. It can be seen that AIDE outperforms both INEXACTDANE and COCOA+. We see similar behavior on other datasets (see Section E.1 of the Appendix). The benefits of AIDE are particularly pronounced in settings with large condition number $\kappa = L/\lambda > N$. This can also be explained theoretically, as convergence rate of fast accelerated stochastic methods has a dependence of $O(\sqrt{N\kappa})$ on $\kappa$ and $N$, as opposed to $O(\kappa)$ without acceleration [18, 10].

In the two left-hand side plots of the bottom row of Figure 1, we demonstrate the effect of local computation effort at each iteration, to solve the local subproblem (2), on the overall convergence of the algorithm. The different lines signify $\{1/6, 1/3, 1, 2, 4\}$ passes through data available locally per iteration of INEXACTDANE. From the plots, it can be seen that running SVRG beyond 4 passes through local data only provides little improvement; thus, suggesting a natural computational setting for INEXACTDANE and demonstrating its advantage over DANE. Finally, in the two right-hand side plots of the bottom row, we show performance of INEXACTDANE with increasing number of nodes. In particular, we partition the data to $\{4, 8, 16, 32, 64\}$ nodes, and keep the number of local iterations of SVRG constant for all settings. The results suggest that in practice, INEXACTDANE can scale gracefully with the number of nodes. We only observe drop in the performance with rcv1 dataset on 64 nodes, where the optimal stepsize of SVRG was significantly smaller than in other cases. Due to space constraints, we relegate experiments on other datasets and loss functions, and demonstration of the robustness of the INEXACTDANE method over DANE to Section E of the Appendix.

## 8  Discussion

In this section, we give a brief comparison of the results we presented. In particular, we compare the key aspects of DANE, DISCO, COCOA+, INEXACTDANE and AIDE algorithms.

---
[2]Datasets available at https://www.csie.ntu.edu.tw/~cjlin/libsvmtools/datasets/

**Communication Complexity**: We showed that AIDE nearly achieves the lower bounds proved in [2] in specific cases. When the functions are quadratic $\delta$-related, AIDE and DISCO match the communication complexity lower bound of $\tilde{O}(\sqrt{\delta/\lambda}\log(1/\epsilon))$. However, DANE and INEXACTDANE have a sub-optimal communication complexity of $\tilde{O}((\delta^2/\lambda^2)\log(1/\epsilon))$ in this case. Similarly, for the general unrelated strongly convex function with condition number $\kappa$, AIDE and DISCO enjoy communication complexity of $\tilde{O}(\sqrt{\kappa}\log(1/\epsilon))$ [2]. Again, this complexity is superior to the convergence rate of $O(\kappa \log(1/\epsilon))$ of COCOA+ , INEXACTDANE and DANE.

**Nature of oracle access**: Both DANE and DISCO require access to a strong oracle. While DANE technically requires an oracle that solves the subproblem (2), DISCO requires a second-order oracle for its execution. On the other hand, both INEXACTDANE and AIDE can be executed using a simple first-order oracle, matching the major advantage of COCOA+.

**Parallelism and Implementation**: One of the appealing aspect of the DANE, INEXACTDANE, AIDE and COCOA+ is the simplicity of implementation and its suitability for large-scale distributed environments. Each iteration of these algorithms is embarrassingly parallelizable since it involves solving a local objective function at each node. The same cannot be said about DISCO due to the asymmetric workload on the master node at each iteration of the inexact damped Newton iteration.

**Distributed SVRG**: As a by-product of our analysis, we obtain partial convergence guarantees of a distributed version of popular SVRG algorithm, that was observed to perform well in practice [7].

In conclusion, AIDE adopts practical advantages of DANE, but at the same time also achieves the optimal communication complexity in [2]. Furthermore, similar to COCOA+, AIDE and INEXACTDANE provide an efficient way of balancing communication and computation complexity; thereby, providing a powerful framework for solving large-scale distributed optimization in a communication-efficient manner.

# Appendix

## A Analysis of INEXACTDANE

## Quadratic case

Before proving the main result, we will need to establish two lemmas.

**Lemma 1.** *Let assumptions of Theorem 1 be satisfied. Then*

$$\left\| w^t - \frac{1}{K} \sum_{k=1}^{K} \hat{w}_k^t \right\| \leq \frac{\eta\gamma}{K} \sum_{k=1}^{K} \|(H_k + \mu I)^{-1} H\| \|w^{t-1} - \hat{w}\|$$

*Proof.* Since $\hat{w}_k^t = \arg\min_w g_{k,\mu}^t(w)$, we have $\nabla g_{k,\mu}^t(\hat{w}_k^t) = 0$. That is, $0 = \nabla g_{k,\mu}^t(\hat{w}_k^t) = \nabla F_k(\hat{w}_k^t) + \eta \nabla f(w^{t-1}) - \nabla F_k(w^{t-1}) + \mu(\hat{w}_k^t - w^{t-1})$. By rearranging the terms in the last expression, we get

$$\hat{w}_k^t = w^{t-1} - \eta \left( (H_k + \mu I)^{-1} H w^{t-1} + (H_k + \mu I)^{-1} l \right). \tag{5}$$

Taking norms on both sides of (5), and using the fact that $\hat{w} = -H^{-1} l$, we obtain

$$\|w^{t-1} - \hat{w}_k^t\| = \eta \|(H_k + \mu I)^{-1} H (w^{t-1} - \hat{w})\|. \tag{6}$$

We proceed to prove the desired result by using the following set of inequalities:

$$\begin{aligned}
\left\| w^t - \frac{1}{K} \sum_{k=1}^{K} \hat{w}_k^t \right\| &\leq \frac{1}{K} \sum_{k=1}^{K} \|w_k^t - \hat{w}_k^t\| \\
&\leq \frac{1}{K} \sum_{k=1}^{K} \gamma \|w^{t-1} - \hat{w}_k^t\| \\
&= \frac{\eta\gamma}{K} \sum_{k=1}^{K} \|(H_k + \mu I)^{-1} H (w^{t-1} - \hat{w})\| \\
&\leq \frac{\eta\gamma}{K} \sum_{k=1}^{K} \|(H_k + \mu I)^{-1} H\| \|w^{t-1} - \hat{w}\|. \tag{7}
\end{aligned}$$

The first inequality follows from Jensen's inequality and the fact that $w^t = \sum_{k=1}^{K} w_k^t / K$. The second inequality follows from the approximation condition on $w_k^t$. The equality is obtained from Equation (6). The last inequality follows from properties of spectral norm of a matrix. □

**Lemma 2.** *Let assumptions of Theorem 1 be satisfied. Then*

$$\|w^t - \hat{w}\| - \left\| \eta \tilde{H}^{-1} H - I \right\| \|w^{t-1} - \hat{w}\| \leq \left\| w^t - \frac{1}{K} \sum_{k=1}^{K} \hat{w}_k^t \right\|.$$

*Proof.*

$$\left\|w^t - \frac{1}{K}\sum_{k=1}^{K}\hat{w}_k^t\right\| = \left\|w^t - w^{t-1} + \frac{\eta}{K}\sum_{k=1}^{K}(H_k+\mu I)^{-1}H(w^{t-1}-\hat{w})\right\|$$
$$= \left\|(w^t-\hat{w}) - (w^{t-1}-\hat{w}) + \eta\tilde{H}^{-1}H(w^{t-1}-\hat{w})\right\|$$
$$= \left\|(w^t-\hat{w}) - (I - \eta\tilde{H}^{-1}H)(w^{t-1}-\hat{w})\right\|$$
$$\geq \|w^t-\hat{w}\| - \left\|(I-\eta\tilde{H}^{-1}H)(w^{t-1}-\hat{w})\right\|$$
$$\geq \|w^t-\hat{w}\| - \left\|I-\eta\tilde{H}^{-1}H\right\|\|w^{t-1}-\hat{w}\|. \quad (8)$$

The first equality follows from the optimality property of $\hat{w}_k^t$ (see Equation (5)). The first inequality follows from the triangle inequality. The second inequality follows from the properties of the matrix spectral norm. □

### A.1 Proof of Theorem 1

*Proof.* Lemmas 1 and 2 imply the inequality

$$\left\|w^t-\hat{w}\right\| - \left\|\eta\tilde{H}^{-1}H - I\right\|\left\|w^{t-1}-\hat{w}\right\| \leq \frac{\eta\gamma}{K}\sum_{k=1}^{K}\left\|(H_k+\mu I)^{-1}H\right\|\left\|w^{t-1}-\hat{w}\right\|.$$

It suffices to rearrange the terms. □

### A.2 Proof of Theorem 2

*Proof.* Follows by using Lemma 4 (see Appendix C) together with Theorem 1. □

### A.3 INEXACTDANE for Stochastic Quadratic Setting

In the interesting case of stochastic quadratic setting (see [20]), we can provide a more precise result. More specifically, we have the following key result in the stochastic quadratic setting.

**Corollary 4.** *Suppose each $H_k$ is the average of $n$ i.i.d. symmetric positive definite matrices with largest eigenvalue bounded above by $L$. Then for any $0 < \alpha \leq 1$, with probability $1-\alpha$ we have:*

i *Functions $\{F_k\}$ are $\delta$-related with $\delta = \sqrt{\frac{32L^2\log(Kd/\alpha)}{n}}$,*

ii *The iterates of INEXACTDANE (Algorithm 1) (with parameters defined in Theorem 2) after*

$$t = \tilde{O}\left(\frac{(L/\lambda)^2}{n}\log\left(\frac{Kd}{\alpha}\right)\log\left(\frac{L\|w^0-\hat{w}\|^2}{\epsilon}\right)\right)$$

*satisfy $f(w^t) \leq f(\hat{w}) + \epsilon$.*

*Proof.* The first part of the proof is directly from Lemma 2 of [20]. The second part of the proof follows from Theorem 2. □

# Strongly Convex Case

## A.4 Proof of Theorem 3

*Proof.* We first observe that

$$\|\nabla g_{k,\mu}^t(w_k^t) - \nabla g_{k,\mu}^t(w^{t-1})\| \geq \|\nabla g_{k,\mu}^t(w^{t-1})\| - \|\nabla g_{k,\mu}^t(w_k^t)\|$$
$$\geq (1-\gamma)\|\nabla g_{k,\mu}^t(w^{t-1})\|. \qquad (9)$$

The first inequality follows from triangle inequality. The second inequality follows from the condition that $\|\nabla g_{k,\mu}^t(w_k^t)\| \leq \gamma \|\nabla g_{k,\mu}^t(w^{t-1})\|$. We now define function $h_k^t : \mathbb{R}^d \to \mathbb{R}$ as

$$h_k^t(w) := F_k(w) + \frac{\mu}{2}\|w - w^{t-1}\|^2.$$

We define the virtual iterate $\hat{w}_k^t = \arg\min_w g_{k,\mu}^t(w)$. Using the optimality conditions of $\hat{w}_k^t$, we have the following:

$$\nabla h_k^t(\hat{w}_k^t) - \nabla h_k^t(w^{t-1}) = -\eta \nabla f(w^{t-1}). \qquad (10)$$

Also, define Bregman divergence of a smooth function $\phi$ as $D_\phi(w, w') = \phi(w) - \phi(w') - \langle \nabla \phi(w'), w - w' \rangle$. We observe the following:

$$f(w_k^t) = f(w^{t-1}) + \langle \nabla f(w^{t-1}), w_k^t - w^{t-1} \rangle + D_f(w_k^t, w^{t-1})$$
$$= f(w^{t-1}) - \frac{1}{\eta}\langle \nabla h_k^t(\hat{w}_k^t) - \nabla h_k^t(w^{t-1}), w_k^t - w^{t-1} \rangle + D_f(w_k^t, w^{t-1})$$
$$= f(w^{t-1}) - \frac{1}{\eta}\underbrace{\langle \nabla h_k^t(\hat{w}_k^t) - \nabla h_k^t(w_k^t), w_k^t - w^{t-1} \rangle}_{T_1}$$
$$- \frac{1}{\eta}\langle \nabla h_k^t(w_k^t) - \nabla h_k^t(w^{t-1}), w_k^t - w^{t-1} \rangle + D_f(w_k^t, w^{t-1}). \qquad (11)$$

The second equality is due to optimality condition in (10). We bound the term $T_1$ in the following manner:

$$|T_1| \leq \|\nabla h_k^t(\hat{w}_k^t) - \nabla h_k^t(w_k^t)\|\|w_k^t - w^{t-1}\|$$
$$\leq (L+\mu)\|\hat{w}_k^t - w_k^t\|\|w_k^t - w^{t-1}\|. \qquad (12)$$

The first inequality is due to Cauchy-Schwarz inequality. The second inequality follows from $L + \mu$ lipschitz continuous nature of the gradient of $h_k^t$. In order to proceed further, a bound on the term $\|w_k^t - w^{t-1}\|$ is obtained in the following fashion:

$$\|w^{t-1} - w_k^t\| \leq \|w^{t-1} - \hat{w}_k^t\| + \|\hat{w}_k^t - w_k^t\|$$
$$\leq \|w^{t-1} - \hat{w}_k^t\| + \gamma\|\hat{w}_k^t - w^{t-1}\| = (1+\gamma)\|w^{t-1} - \hat{w}_k^t\|. \qquad (13)$$

The first inequality is due to triangle inequality. The second inequality is due to the inexactness condition in Algorithm 1. Substituting the above bound in Equation (12), we have the following:

$$|T_1| \leq (L+\mu)(1+\gamma)\|\hat{w}_k^t - w_k^t\|\|w^{t-1} - \hat{w}_k^t\| \leq (L+\mu)(1+\gamma)\gamma\|w^{t-1} - \hat{w}_k^t\|^2. \qquad (14)$$

The second inequality is again due to inexactness condition in Algorithm 1. In order to bound $\|w^{t-1} - \hat{w}_k^t\|^2$, we observe the following:

$$(\mu + \lambda)\|w^{t-1} - \hat{w}_k^t\| \leq \|\nabla g_{k,\mu}^t(w^{t-1})\| = \|\eta \nabla f(w^{t-1})\|. \qquad (15)$$

The first inequality follows from Lemma 5. Using this relationship in Equation (14), we have the following:

$$|T_1| \leq \frac{2\gamma\eta^2(L+\mu)}{(\lambda+\mu)^2}\|\nabla f(w^{t-1})\|^2.$$

Substituting the bound in Equation (11), we have:

$$f(w_k^t) \leq f(w^{t-1}) + \frac{2\gamma\eta(L+\mu)}{(\lambda+\mu)^2}\|\nabla f(w^{t-1})\|^2$$
$$- \frac{1}{\eta(L+\mu)}\|\nabla h_k^t(w_k^t) - \nabla h_k^t(w^{t-1})\|^2 + D_f(w_k^t, w^{t-1})$$
$$\leq f(w^{t-1}) + \frac{2\gamma\eta(L+\mu)}{(\lambda+\mu)^2}\|\nabla f(w^{t-1})\|^2$$
$$- \frac{1}{\eta(L+\mu)}\|\nabla h_k^t(w_k^t) - \nabla h_k^t(w^{t-1})\|^2 + \frac{L}{2}\|w^{t-1} - w_k^t\|^2$$
$$\leq f(w^{t-1}) + \frac{2\gamma\eta(L+\mu)}{(\lambda+\mu)^2}\|\nabla f(w^{t-1})\|^2$$
$$- \frac{1}{\eta(L+\mu)}\|\nabla g_{k,\mu}^t(w_k^t) - \nabla g_{k,\mu}^t(w^{t-1})\|^2 + \frac{2L\eta^2}{(\lambda+\mu)^2}\|\nabla f(w^{t-1})\|^2$$
$$\leq f(w^{t-1}) - \left[\frac{(1-\gamma)^2}{\eta(L+\mu)} - \frac{2L}{(\lambda+\mu)^2} - \frac{2\gamma(L+\mu)}{\eta(\lambda+\mu)^2}\right]\eta^2\|\nabla f(w^{t-1})\|^2$$
$$\leq f(w^{t-1}) - \left[\frac{(1-\gamma)^2}{\eta(L+\mu)} - \frac{2L}{(\lambda+\mu)^2} - \frac{2\gamma(L+\mu)}{\eta(\lambda+\mu)^2}\right]\eta^2\lambda(f(w^{t-1}) - f(\hat{w})). \tag{16}$$

The first step is due to Lemma 5. The second step follows the $L$-smoothness of $f$. The third step follows from the fact that $\|\nabla g_{k,\mu}^t(w_k^t) - \nabla g_{k,\mu}^t(w^{t-1})\| = \|\nabla h_k^t(w_k^t) - \nabla h_k^t(w^{t-1})\|$ and bound on $\|w^{t-1} - w_k^t\|^2$ from Equation (13) and (15). The fourth inequality follows from Equation (9) and the fact that $\nabla g_{k,\mu}^t(w^{t-1}) = \nabla f(w^{t-1})$. The last inequality is due to strong convexity of function $f$. Finally, we observe:

$$(f(w^t) - f(\hat{w}) \leq \frac{1}{K}\sum_{k=1}^{K}(f(w_k^t) - f(\hat{w})) \leq (1-\tilde{\rho})(f(w^{t-1}) - f(\hat{w})).$$

The first inequality is due to convexity. The second inequality follows from Equation (16). Therefore, we have the required result. $\square$

## Weakly Convex Case

### Proof of Corollary 2

We observe the following:

$$f(w^t) \leq f_\epsilon(w^t) \leq \min_w f_\epsilon(w) + (1-\rho_\epsilon)^s[f_\epsilon(w^0) - \min_w f_\epsilon(w)]$$
$$\leq f(\hat{w}) + \frac{\epsilon}{2}\|\hat{w} - w^0\|^2 + (1-\rho_\epsilon)^s[f_\epsilon(w^0) - \min_w f_\epsilon(w)]$$
$$\leq f(\hat{w}) + \epsilon\left[(1/2)\|\hat{w} - w^0\|^2 + (f(w^0) - f(\hat{w}))\right]. \tag{17}$$

The second inequality is a consequence of Equation (4). The last inequality follows from the facts that (i) $s = \tilde{O}(L\log(1/\epsilon)/\epsilon)$ and (ii) $[f_\epsilon(w^0) - \min_w f_\epsilon(w)] = [f(w^0) - \min_w f_\epsilon(w)] \leq [f(w^0) - f(\hat{w})]$.

# Nonconvex Case

## A.5 Proof of Theorem 4

*Proof.* Using essentially the same argument as in Theorem 3, we have the following:

$$f(w_k^t) \leq f(w^{t-1}) - \left[\frac{(1-\gamma)^2}{\eta(L+\mu)} - \frac{2L}{(\mu-L)^2} - \frac{2\gamma(L+\mu)}{\eta(\mu-L)^2}\right] \eta^2 \|\nabla f(w^{t-1})\|^2.$$

The argument uses the fact that $g_{k,\mu}^t$ is $(\mu+L)$-smooth and $(\mu-L)$-strongly convex. This is in turn due to the lipschitz continuous nature of the gradient. Note that this is possible only when $\mu > L$ (as mentioned in the theorem statement). Using the definition of $w^t$ in the nonconvex case, we have

$$f(w^t) \leq f(w^{t-1}) - \left[\frac{(1-\gamma)^2}{\eta(L+\mu)} - \frac{2L}{(\mu-L)^2} - \frac{2\gamma(L+\mu)}{\eta(\mu-L)^2}\right] \eta^2 \|\nabla f(w^{t-1})\|^2. \tag{18}$$

Rearranging the above inequality and using the definition of $\theta$, we get

$$\|\nabla f(w^{t-1})\|^2 \leq \frac{f(w^{t-1}) - f(w^t)}{\theta}.$$

Using a telescopic sum on the above inequality and the fact that $\hat{w}$ is an optimal solution of (1), we have the desired result. $\square$

# B Analysis of AIDE

# Quadratic Case

## B.1 Proof of Theorem 5

*Proof.* First consider the case where $2\sqrt{2}\delta \leq \lambda$. In this case, we have $\tau = 0$ and $1/16 \leq \gamma \leq 1/8$. The required result trivially follows from Theorem 2. We turn our attention to the interesting case of $2\sqrt{2}\delta > \lambda$. For this case, we choose $\tau = 2\sqrt{2}\delta - \lambda$. We use $F_k^t(w)$ to denote the function $F_k(w) + (\tau/2)\|w - y^{t-1}\|^2$. As mentioned in pseudocode of Algorithm 2, we use INEXACTDANE for the set of functions $F_k^t$. Let $\hat{u}^t = \arg\min_w f^t(w)$ (refer to Algorithm 2). By solving each subproblem inexactly for $s$ iterations, from Theorem 2, we have the following:

$$f^t(w^t) - f^t(\hat{u}^t) \leq \tfrac{L+\tau}{2} \|w^t - \hat{u}^t\|^2 \leq \tfrac{L+\tau}{2} \tilde{\rho}^s \|w^{t-1} - \hat{u}^t\|^2, \tag{19}$$

where $\tilde{\rho} = 2/3$. The first inequality follows from the $L + \tau$ Lipschitz continuous gradient of $f^t$. The second inequality follows from Theorem 2. The value of $\tilde{\rho}$ is obtained from Theorem 2 noting the following two facts. (i) The inexactness parameter $\gamma \leq 1/8$. (ii) $\tilde{\rho} = 2/3$ when $\gamma \leq 1/8$ and $2\sqrt{2}\delta \leq (\lambda + \tau)$. Here, we used the fact that INEXACTDANE is applied on the function $f^t$ and $\nabla^2 f^t \succeq (\lambda + \tau)I$. Using Equation (19), we have the following:

$$f^t(w^t) - f^t(\hat{u}^t) \leq \tfrac{L+\tau}{2} \tilde{\rho}^s \|w^{t-1} - \hat{u}^t\|^2 \leq \tfrac{L+\tau}{\lambda+\tau} \tilde{\rho}^s [f^t(w^{t-1}) - f^t(\hat{u}^t)].$$

This simply follows from the fact that $\nabla^2 f^t \succeq (\lambda + \tau)I$. Using Proposition 3.2 of [10], we have the total number of iterations of INEXACTDANE iterations one has to execute to achieve $\epsilon$-accuracy is

$$t = \tilde{O}\left(\tfrac{1}{1-\tilde{\rho}} \sqrt{\tfrac{\lambda+\tau}{\lambda}} \log(1/\epsilon)\right) = \tilde{O}\left(\sqrt{\tfrac{\delta}{\lambda}} \log(1/\epsilon)\right).$$

This is obtained from the fact that $\tilde{\rho} = 2/3$. Therefore, we get the desired result. $\square$

## B.2 AIDE in the Stochastic Quadratic Setting

The following result follows as a Corollary of Theorem 5.

**Corollary 5.** *Suppose each $H_k$ is the average of $n$ i.i.d positive semidefinite matrices with eigenvalues at most $L$. Then with probability at least $1 - \alpha$:*

*i Functions $\{F_k\}$ are $\delta$-related with $\delta = \sqrt{\frac{32L^2 \log(Kd/\alpha)}{n}}$,*

*ii The iterates of INEXACTDANE (Algorithm 1) (with parameters in Theorem 5) after*

$$t = \tilde{O}\left(\frac{\sqrt{L}}{n^{1/4}\sqrt{\lambda}} \log^{1/4}\left(\frac{Kd}{\alpha}\right) \log\left(\frac{L\|w^0 - \hat{w}\|^2}{\epsilon}\right)\right)$$

*satisfy $f(w^t) \leq f(\hat{w}) + \epsilon$.*

*Proof.* The first part of the proof is directly from Lemma 2 of [20]. The second part of the proof follows from Theorem 5. □

# Strongly convex case

## B.3 Proof of Theorem 6

*Proof.* Let $\hat{u}^t = \arg\min_w f^t(w)$. We first observe that after $s$ iterations of INEXACTDANE in the $t^{\text{th}}$ iteration of AIDE (applied on $f^t$), we have the following:

$$f^t(w^t) - f^t(\hat{u}^t) \leq (1 - \tilde{\rho})^s (f^t(w^{t-1}) - f^t(\hat{u}^t)), \tag{20}$$

where

$$\tilde{\rho} = \left[\frac{49}{64(L + \mu + \tau)} - \frac{2(L + \tau)}{(\lambda + \tau + \mu)^2} - \frac{(L + \tau + \mu)}{4(\lambda + \tau + \mu)^2}\right] L$$

This follows directly from Theorem 3 and the fact that $f^t$ has $(L + \tau)$-Lipschitz continuous gradient and is $(\lambda + \tau)$-strongly convex. Note that with the given value of $\tau$, $\mu$ and $\gamma$, we have following:

$$\tilde{\rho} = L\left[\frac{49}{64(2L + \mu - \lambda)} - \frac{2(2L - \lambda)}{(L + \mu)^2} - \frac{2L + \mu - \lambda}{4(L + \mu)^2}\right] \tag{21}$$

$$\geq L\left[\frac{49}{64(14L)} - \frac{4L}{(13L)^2} - \frac{14L}{4(13L)^2}\right] \geq 0.01. \tag{22}$$

The above bound follows from simple algebra. Now again using Proposition 3.2 of [10], the total number of INEXACTDANE iterations required to achieve an $\epsilon$-accurate solution is

$$t = \tilde{O}\left(\frac{1}{\tilde{\rho}}\sqrt{\frac{\lambda + \tau}{\lambda}} \log\left(\frac{1}{\epsilon}\right)\right) = \tilde{O}\left(\sqrt{\frac{L}{\lambda}} \log\left(\frac{1}{\epsilon}\right)\right).$$

This is obtained from the bound in Equation (21). This gives us the desired result. Note that the constants in this result are not optimized. □

# C Auxiliary Results

Here we establish two lemmas which are used in the proofs of the main results.

The following result is a slight extension of Lemma 4 in [20].

**Lemma 3.** *Let $\mu, \delta, \xi$ be positive constants satisfying $\max\{\mu, \delta\} < \xi$. Further, let $A$ and $B_1, \ldots, B_K$ be $d \times d$ symmetric matrices satisfying $\xi I \preceq A$, $\sum_k B_k = 0$ and $\|B_k\| \leq \delta$ for all $k \in [K]$. Then, we have the following:*

$$\left\|\left((A+B_k)^{-1} - A^{-1}\right)(A - \mu I)\right\| \leq \frac{2\delta}{\xi - \delta}, \quad k \in [K],$$

$$\left\|\left(\frac{1}{K}\sum_{k=1}^K (A+B_k)^{-1} - A^{-1}\right)(A - \mu I)\right\| \leq \frac{2\delta^2}{\xi(\xi - \delta)}.$$

*Proof.* First, we observe that

$$(A + B_k)^{-1} = (A(I + A^{-1}B_k))^{-1} = (I + A^{-1}B_k)^{-1}A^{-1} = A^{-1} + \sum_{j=1}^{\infty}(-1)^j (A^{-1}B_k)^j A^{-1}.$$

The above equality follows form series expansion of $(I + A^{-1}B_k)^{-1}$, which is possible as $\|A^{-1}B_k\| \leq \|A^{-1}\|\|B_k\| = (\delta/\xi) < 1$. From the above equality, we obtain the following bound:

$$\left\|\left((A+B_k)^{-1} - A^{-1}\right)(A - \mu I)\right\| = \left\|\sum_{j=1}^{\infty}(-1)^j(A^{-1}B_k)^j A^{-1}(A - \mu I)\right\|$$

$$\leq \sum_{j=1}^{\infty}\left\|(-1)^j(A^{-1}B_k)^j A^{-1}(A - \mu I)\right\|$$

$$\leq \sum_{j=1}^{\infty}\left\|A^{-1}\right\|^j \|B_k\|^j \left\|I - \mu A^{-1}\right\|$$

$$\leq \sum_{j=1}^{\infty} \frac{\delta^j}{\xi^j}\left(1 + \frac{\mu}{\xi}\right) \leq \frac{2\delta}{\xi - \delta}.$$

The first inequality follows from the triangle inequality. The second inequality follows from the Cauchy-Schwarz inequality. The last inequality follows from the fact that $\mu/\xi < 1$.

The second part of the claim was established as Lemma 4 in [20]. □

**Lemma 4.** *Assume that $\lambda I \preceq H \preceq LI$, where $H := \frac{1}{K}\sum_{k=1}^K H_k$ and $\lambda > 0$. Further, assume that $\|H_k - H\| \leq \delta$ for all $k \in [K]$ and let $\tilde{H}^{-1} := \frac{1}{K}\sum_{k=1}^K (H_k + \mu I)^{-1}$, where $\mu = \max\left\{0, 8\frac{\delta^2}{\lambda} - \lambda\right\}$. If we let*

$$\gamma := \begin{cases} \frac{1}{8} & \text{if } 2\sqrt{2}\delta \leq \lambda, \\ \frac{\lambda^2}{192\delta^2} & \text{otherwise,} \end{cases}$$

*then*

$$\rho = \|\tilde{H}^{-1}H - I\| + \frac{\gamma}{K}\sum_{k=1}^K \|(H_k + \mu I)^{-1}H\| \leq \begin{cases} \frac{2}{3} & \text{if } 2\sqrt{2}\delta \leq \lambda, \\ 1 - \frac{\lambda^2}{24\delta^2} & \text{otherwise.} \end{cases}$$

*Proof.* Using Lemma 3 with $A = H + \mu I$ and $B_k = H_k - H$, we have the following inequality:

$$\left\|\left((H_k + \mu I)^{-1} - (H + \mu I)^{-1}\right)H\right\| \leq \frac{2\delta}{\lambda + \mu - \delta}, \tag{23}$$

$$\left\|\left(\tilde{H}^{-1} - (H + \mu I)^{-1}\right)H\right\| \leq \frac{2\delta^2}{(\lambda + \mu)(\lambda + \mu - \delta)}, \tag{24}$$

for all $k \in [K]$. From the above inequalities, we have

$$\begin{aligned}
\|(H_k + \mu I)^{-1} H\| &\leq \|(H_k + \mu I)^{-1} H - I\| + \|I\| \\
&\leq \|(H_k + \mu I)^{-1} H - (H + \mu I)^{-1} H\| + \|(H + \mu I)^{-1} H - I\| + \|I\| \\
&\leq \frac{2\delta}{\lambda + \mu - \delta} + \frac{\mu}{\lambda + \mu} + 1.
\end{aligned} \quad (25)$$

The first and second inequality follow from triangle inequality. The third inequality follows from (23) and Lemma 3 in [20], which says that $\|(H + \mu I)^{-1} H - I\| = \mu/(\lambda + \mu)$. Furthermore, we also have the following bound:

$$\begin{aligned}
\|\tilde{H}^{-1} H - I\| &\leq \|\tilde{H}^{-1} H - (H + \mu I)^{-1} H\| + \|(H + \mu I)^{-1} H - I\| \\
&\leq \frac{2\delta^2}{(\lambda + \mu)(\lambda + \mu - \delta)} + \frac{\mu}{\lambda + \mu}.
\end{aligned} \quad (26)$$

The first and second inequality follow from triangle inequality and Equation (24), respectively. Inequality in Equation (26) was earlier shown in [20]. Using inequalities in Equation (25) and Equation (26), we obtain the following bound:

$$\begin{aligned}
\rho &= \|\tilde{H}^{-1} H - I\| + \frac{\gamma}{K} \sum_{k=1}^{K} \|(H_k + \mu I)^{-1} H\| \\
&\leq \left[ \frac{2\delta^2}{(\lambda + \mu)(\lambda + \mu - \delta)} + \frac{\mu}{\lambda + \mu} \right] + \gamma \left[ \frac{2\delta}{\lambda + \mu - \delta} + \frac{\mu}{\lambda + \mu} + 1 \right]
\end{aligned}$$

Let us first consider the case of $\lambda \geq 2\sqrt{2}\delta$. In this case, we set $\mu = 0$, $\gamma = 1/8$ and hence, we have

$$\rho \leq \frac{2\delta^2}{\lambda(\lambda - \delta)} + \gamma \left[ \frac{2\delta}{\lambda - \delta} + 1 \right] < 0.53 < \frac{2}{3}.$$

Now, consider the case of $\lambda < 2\sqrt{2}\delta$. We use $\psi$ to denote $8\delta/\lambda$. Note that $\psi > 2$. In this case, we set $\mu = (8\delta^2/\lambda) - \lambda$ and $\gamma = \lambda^2/(192\delta^2)$ and hence, the following holds:

$$\rho \leq \frac{2}{\psi(\psi - 1)} + (1 + \gamma)\left(1 - \frac{\lambda}{\psi \delta}\right) + \gamma\left(1 + \frac{2}{\psi - 1}\right). \quad (27)$$

We observe that $\psi - 1 > \psi/2$ (since $\psi > 2$). Further, the following inequality holds:

$$(1 + \gamma)\left(1 - \frac{\lambda}{\psi \delta}\right) \leq 1 - \frac{\lambda}{\psi \delta} + \gamma.$$

Substituting the above inequalities in Equation (27), we have the following:

$$\rho \leq \frac{4}{\psi^2} + 1 - \frac{\lambda}{\psi \delta} + \gamma + \gamma\left(1 + \frac{2}{\psi - 1}\right) \leq 1 - \frac{\lambda^2}{16\delta^2} + 4\gamma \leq 1 - \frac{\lambda^2}{24\delta^2}.$$

The second inequality is due to that fact that $\psi - 1 > 1$ and the fact that $\gamma = \lambda^2/(192\delta^2)$. Hence, we have the desired result. $\square$

**Lemma 5.** *Suppose a function $h : \mathbb{R}^d \to \mathbb{R}$ is $L$-smooth and $\lambda$-strongly convex. Let $w^* = \arg\min h(w)$. Then, we have the following:*

$$\langle \nabla h(w') - h(w), w' - w \rangle \geq \frac{1}{L} \|\nabla h(w') - \nabla h_i(w)\|^2,$$
$$\|\nabla h(w)\| \geq \lambda \|w - w^*\|,$$

*for all $w', w \in \mathbb{R}^d$.*

**Algorithm 3:** A DISTRIBUTED VERSION OF SVRG$(f, w^0, s, r, h)$
---
**Data:** $f(w) = \frac{1}{K} \sum_{k=1}^{K} F_k(w)$, initial point $w^0 \in \mathbb{R}^d$, stepsize $h > 0$
**for** $t = 1$ *to* $s$ **do**
    Compute $\nabla f(w^t)$ and distribute to all machines
    **for** $k = 1$ *to* $K$ **do in parallel**
        $w_k^t = w^t$
        **for** $j = 1$ *to* $r$ **do**
            Sample $i \in \mathcal{P}_k$ (e.g., uniformly at random)
            $w_k^t = w_k^t - h \left( \nabla f_i(w_k^t) - \nabla f_i(w^{t-1}) + \nabla f(w^{t-1}) \right)$
        **end**
    **end**
    $w^t = \frac{1}{K} \sum_{k=1}^{K} w_k^t$
**end**
**return** $w^t$
---

## D  A Practical Distributed Version of SVRG

In Section 6 we pointed out how running SVRG as a local solver in INEXACTDANE in certain setting is equivalent to the Algorithm 3 below. It is a distributed version of SVRG, that has been observed to perform well in practice [7], but has not been directly analyzed. Note that there exist another way to obtain exactly the same algorithm. That is to rewrite the local subproblem (2) as follows

$$g_{k,\mu}^t(w) = \sum_{i \in \mathcal{P}_k} \left[ f_i(w) - \langle \nabla f_i(w^{t-1}) - \nabla f(w^{t-1}), w \rangle + \frac{\mu}{2} \| w - w^{t-1} \|^2 \right],$$

and directly apply SGD locally on this reformulation within INEXACTDANE and set $\mu = 0$.

## E  Additional Experiments

In this section, we provide extended version of what was already presented in Section 7, along with results of different types. In particular, we extend the comparison of INEXACTDANE, AIDE and COCOA+ to various settings and with different datasets. We also study of the effect of varying local iterations between rounds of communication on the performance of the algorithms. We present further results showing performance under varying number of nodes onto which the dataset is distributed, and highlight a case of non-random data distribution, in which DANE diverges, while the performance of INEXACTDANE degrades only slightly, compared to benchmark with random data distribution.

Again, the following default statements are true, unless stated otherwise. We use SVRG as local solver for INEXACTDANE and AIDE, and we use SDCA in COCOA+. We run all presented methods for a single pass over data available locally in every iteration. We partition data randomly (mimicking the iid data distribution) across 8 nodes. We present the best performance selected from a range of stepsize parameters for SVRG and aggregation parameters for COCOA+. In the plots, we use DANE as label for INEXACTDANE.

We omit the experiments for quadratic objectives as the results were very similar to the ones presented here and hence, did not provide any additional insights, compared to the experiments with classification on publicly available datasets.

### E.1  INEXACTDANE, AIDE and COCOA+

In Figures 2 and 3, we present a comparison of INEXACTDANE, AIDE and COCOA+ on a binary classification task on the rcv1 dataset. Plots in top row are for logistic loss, and bottom is for smoothed hinge loss. We apply L2-regularization with $\lambda$ set to $1/(cN)$ for $c \in \{1, 10, 100\}$, which correspond to left, middle and right column respectively. In Figure 2 we partition the data randomly to 8 nodes, while in Figure 3 we use partitioning across

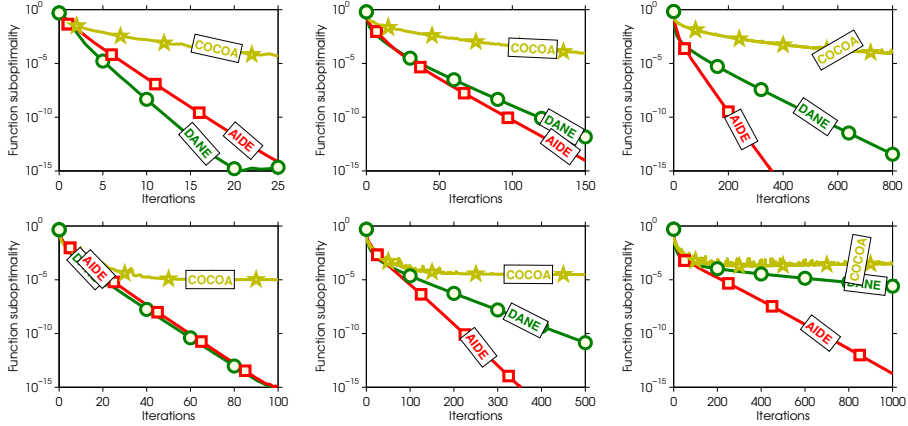

Figure 2: rcv1 dataset, 8 nodes, regularization parameter $\lambda$ set to $1/(cN)$, for $c \in \{1, 10, 100\}$ (in left/middle/right columns respectively). Top row: logistic loss. Bottom row: smoothed hinge loss.

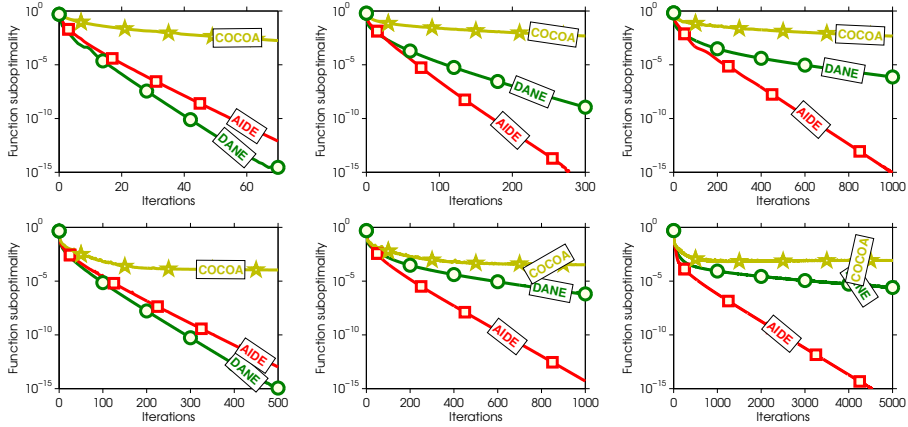

Figure 3: rcv1 dataset, 64 nodes, regularization parameter $\lambda$ set to $1/(cN)$, for $c \in \{1, 10, 100\}$ (in left/middle/right columns respectively). Top row: logistic loss. Bottom row: smoothed hinge loss.

64 nodes. To strengthen our claims, in Figures 4 and 5, we provide experiment with settings analogous to the ones in Figure 2, but on covtype and realsim datasets respectively.

In this experiment, we can see a common pattern arising in different settings. The benefit of acceleration in AIDE is present only when the condition number $\kappa := L/\mu > N$, and the larger it is, the larger is the gap in performance. This is to be expected, as the acceleration of the fast stochastic methods changes $\kappa$ in convergence rates to $\sqrt{N\kappa}$ [18]. Both INEXACTDANE and AIDE outperform COCOA+, with AIDE performing much better in all studied settings. The behaviour of COCOA+, where if one looks only at suboptimality of primal function value, the algorithm quickly converges to modest accuracy, and then converges very slowly to higher accuracies, has been confirmed as correct by its authors.

## E.2 Varying amount of local computation

In Figure 6 we demonstrate how spending more local computational resources in each iteration to solve the local subproblem (2) leads to faster convergence in terms of number of iterations. Note that in the settings presented, it is not possible to get close form solution to the subproblems, and hence we can only approximate behaviour of DANE by running a local method for a long time. In number of cases in the above, running

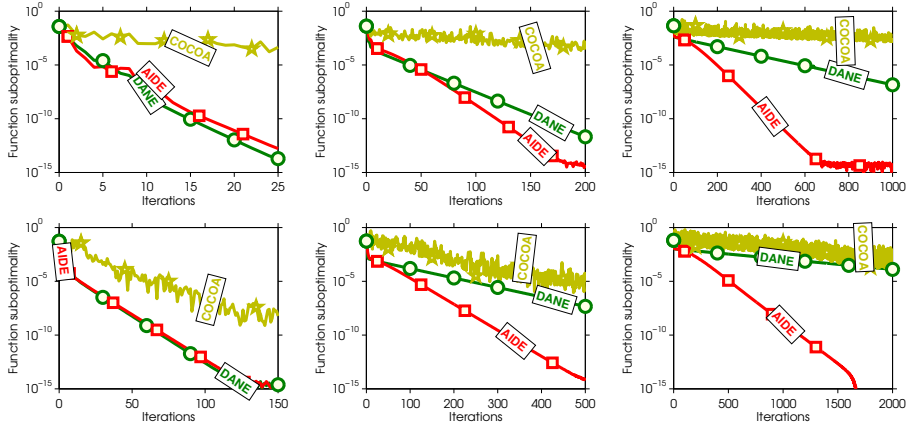

Figure 4: covtype dataset, 8 nodes, regularization parameter $\lambda$ set to $1/(cN)$, for $c \in \{1, 10, 100\}$ (in left/middle/right columns respectively). Top row: logistic loss. Bottom row: smoothed hinge loss.

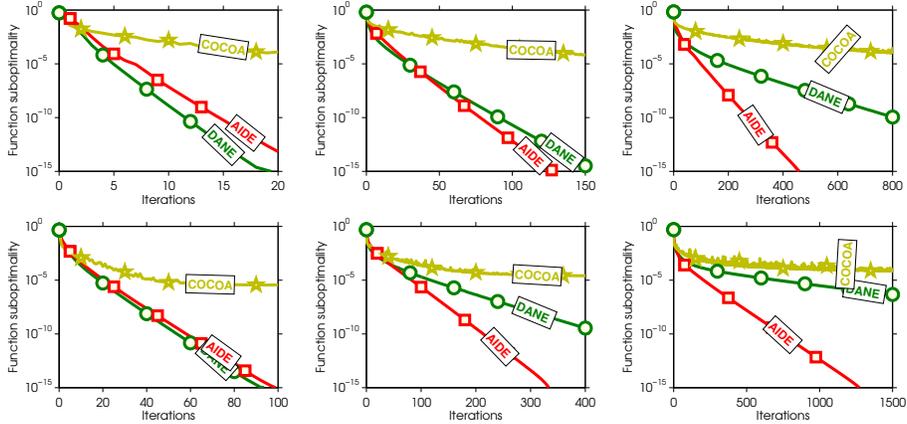

Figure 5: realsim dataset, 8 nodes, regularization parameter $\lambda$ set to $1/(cN)$, for $c \in \{1, 10, 100\}$ (in left/middle/right columns respectively). Top row: logistic loss. Bottom row: smoothed hinge loss.

SVRG for 4 passes through local data already provides little to none improvement. The only dataset on which significant improvement is visible is covtype. This behaviour likely is due to the fact that $N \gg d$, and hence under random data distribution, the local problems can be seen as $\delta$-related with very small $\delta$.

The labels in the figure mean represent the following: The labels DANE/6 and DANE/3 correspond to running the local SVRGalgorithm for one-sixth and one-third of pass through local data in every iteration. The labels DANE*2 and DANE*4 correspond to running the local SVRG algorithm for two and four passes through the local data.

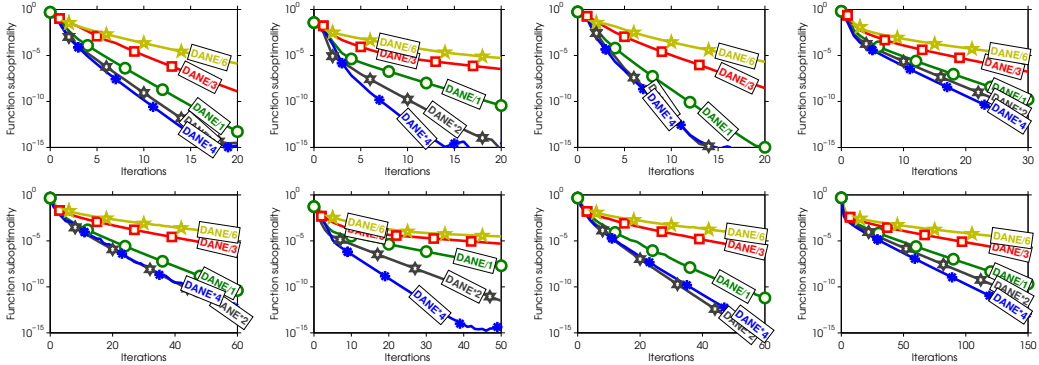

Figure 6: Varying number of passes through local data per iteration in range $[1/6, 4]$. Top row: logistic loss, Bottom row: smoothed hinge loss. Datasets in columns: rcv1, covtype, realsim, url

### E.3 Node scaling

The plots in top row in Figure 7 show how the performance changes, as we partition the data across different number of nodes, but keep the local work equal to one pass through local data. The performance degrades as we increase number of nodes, because we do less relative work on any single computer. However, this is a positive result, as the algorithm *always* converges.

In the bottom row, we double the number of passes over local data when we double the number of nodes. This is motivated to compare how the algorithm performs, when equal number of iteration of local SVRG is used. In most cases, the algorithms perform very similarly, demonstrating the robustness of INEXACTDANE. For rcv1 and realsim datasets partitioned on 64 nodes, the convergence slows down. In this case, the optimal stepsize for local SVRG was significantly smaller than in the other cases. This was likely caused by getting into region where the aggregation starts to be unstable, and thus higher accuracy on local subproblems leads to worse overall performance.

### E.4 Inconvenient data partitioning

In the last experiment, we depart from theory, and observe that INEXACTDANE is, compared to DANE, much more robust to arbitrary partitions of the data. In particular, in Figure 8 we compare INEXACTDANE in two settings: *random*, which corresponds to random data partition to 2 nodes, and *output*, where the data are partitioned according to their output label — positive examples on one node, negative examples on the other. In this setting, DANE diverges, while the performance of INEXACTDANE drops down only slightly. We observed that the optimal stepsize for local SVRG is about $5 - 10\%$ smaller in the *output* case, compared to the *random* case.

This observation is particularly appealing for practitioners, as the data in huge scale applications are often partitioned "as is", i.e., it is given and one does not have the opportunity to randomly reshuffle the data between nodes.

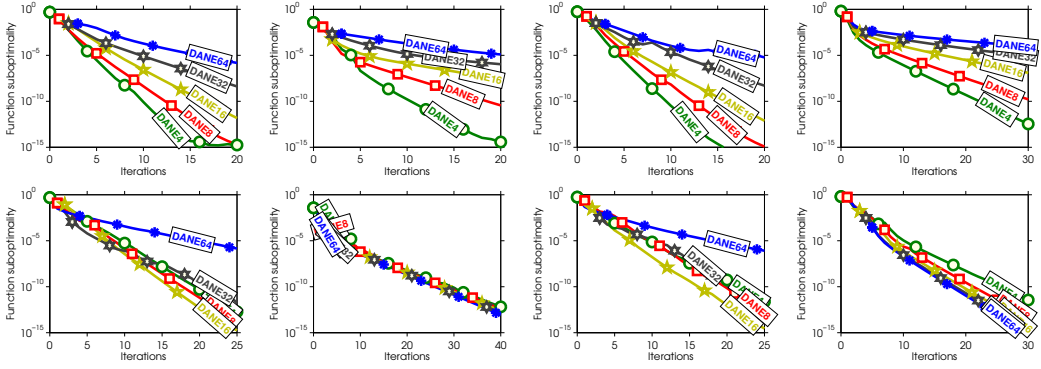

Figure 7: Performance comparison as data is partitioned across 4–64 nodes, with logistic loss. Top row: Single pass through local data per iteration. Bottom row: Fixed number of updates of local SVRG per iteration. Datasets in columns: rcv1, covtype, realsim, url

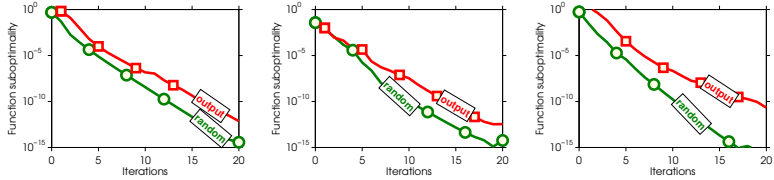

Figure 8: Data partitioned randomly, vs. data partitioned according to their output label. Datasets in columns: rcv1, covtype and realsim.

**Practical considerations:** Although the experiments in Section E.1 suggest superiority of AIDE, it may not always be the case in practice. AIDE comes with the additional requirement to set the catalyst acceleration parameter $\tau$. At the moment, there is not any simple rule guiding its choice, as the optimal $\tau$ depends on properties of the algorithm being accelerated, which are usually not known — see Section 3.1 of [10] for further details. In contrast, COCOA+ is usually a slightly easier to tune in practice, since with SDCA one can use data-independent aggregation parameter equal to $1/K$, and effectively have hyper-parameter-free algorithm. In the case of INEXACTDANE, natural local solvers would require a choice of hyper-parameter such as stepsize. While this can often be affordable to compute, it is not data-independent.

23